\documentclass[12pt]{article}
\usepackage{amsfonts}
\usepackage{amsmath}
\usepackage{nccmath}

\begin{document}

\begin{center}

{\huge Frobenius manifolds associated to \\
\vspace{0,2cm}
Coxeter groups of type $E_7$ and $E_8$} \\
\vspace{0,5cm}
{\large Devis Abriani} \\
\vspace{0,5cm}
SISSA, v.Beirut 2-4, 34151 Trieste, Italy \\
\vspace{0,1cm}
email: abriani@sissa.it \\
\vspace{0,5cm}
{\large October 28, 2009}
\vspace{1cm}

\end{center}

\begin{abstract}
\noindent
Flat coordinates for Frobenius manifolds defined on the orbit space of a Coxeter group $W$ are specified through a certain system of generators of $W$-invariant polynomials. In this note, starting from basic invariants proposed by M.Mehta, we calculate flat coordinates for the exceptional groups of type $E_7$ and $E_8$, leading to a derivation of the potentials for the associated Frobenius structures.
\end{abstract}

\section{Introduction}
A Coxeter group is a group $W$ of linear transformations of a Euclidean space generated by reflections. Irreducible finite Coxeter groups are labelled by $A_l$, $B_l$, $D_l$, $E_6$, $E_7$, $E_8$, $F_4$, $G_2$, $H_3$, $H_4$ and $I_2(p)$ \cite{Cox1}. The orbit space of such groups can be endowed with a structure of Frobenius manifold, due to B.Dubrovin \cite{Du96} and to the previous work by K.Saito, T.Yano and J.Sekiguchi \cite{Sa79,Sa93,SYS}. Flat coordinates, a key ingredient for the structure, correspond to a particular choice of generators of the ring of $W$-invariant polynomials. They have been calculated explicitely by Saito et al. for all irreducible finite Coxeter groups but $E_7$ and $E_8$. In the framework of singularity theory, flat coordinates on the base space of universal unfoldings of isolated hypersurface singularities have been calculated also for $E_7$ and $E_8$ type \cite{Ya,Kat}. The potentials for $E_7$ and $E_8$ type Coxeter groups have been calculated in \cite{DLZ}, but without relating them to an explicit construction of the flat generators. In this note, after briefly recalling the main ingredients of Frobenius manifold structures on the orbit space of finite Coxeter groups, I consider the two cases not treated in \cite{SYS}, deriving explicit expressions for the Saito flat coordinates via generators of the ring of $E_7$- and $E_8$-invariant polynomials proposed by M.Mehta \cite{Meh}. I am very grateful to my supervisor, Prof. B.Dubrovin, for pointing out to me this problem, for fruitful discussions about the subject and for his continuous encouragement and support.

\section{Frobenius manifolds and Coxeter groups}
A Frobenius algebra $(A,\star)$ over a ring $R$ is an associative $R$-algebra with a unity and a symmetric non-degenerate $R$-bilinear inner product $(,)$ such that, for any $a,b,c \in A$:
$$(a\star b,c)=(a,b\star c)$$
A (complex) Frobenius manifold $M$ is a $\mathbb{C}$-manifold with a commutative Frobenius algebra structure over each tangent plane, analytically depending on the coordinates of the point and satisfying a few integrability conditions:
\begin{flushleft}
\begin{enumerate}
\item The metric $\eta_{\alpha\beta}$ on $M$ defined by the inner product $(,)$ is flat;
\item The unity vector field $e$ defined by the algebras is covariantly constant with respect to the Levi-Civita connection defined by $\eta_{\alpha\beta}$ $(\nabla e=0)$;
\item The tensor $\nabla_z (u \star v,w)$ is totally symmetric for any vector fields $u,v,w,z$.
\footnote{The definition of Frobenius manifold usually requires also the existence of a global linear vector field $E$ acting by conformal transformations on the metric and by rescalings on the Frobenius algebras $T_x M$.}
\end{enumerate}
\end{flushleft}
The metric $\eta_{\alpha\beta}$ also defines a contravariant metric $\eta^{\alpha\beta}$ on the cotangent bundle $T_*M$, where another flat metric $g^{\alpha\beta}$ can be constructed \cite{Du96}, giving rise to a flat pencil of metrics $g^{\alpha\beta}+\lambda \eta^{\alpha\beta}$, $\lambda\in\mathbb{C}$. \\
Locally a Frobenius structure can be described, in a suitable set of coordinates $t=(t_1,t_2,\ldots,t_n)$, named flat, by the third derivatives of a function on the manifold $F(t)$ called Frobenius potential.
In such coordinates the metric $\eta_{\alpha\beta}$ is constant and is given by:
\begin{displaymath}
\eta_{\alpha\beta}(t)=\frac{\partial^3 F(t)}{\partial t_1\partial t_{\alpha}\partial t_{\beta}},
\end{displaymath}
where $t_1$ is such that $e=\frac{\partial}{\partial t_1}$, while the structure constants of the Frobenius algebras on $T_t M$ are given by:
\begin{displaymath}
c^{\gamma}_{\alpha\beta}(t)=\eta^{\gamma\delta}(t)\frac{\partial^3 F(t)}{\partial t_{\delta} \partial t_{\alpha}\partial t_{\beta}}.
\end{displaymath}
\indent
A Coxeter group $W$ is a group of linear transformations on a Euclidean space $V$ generated by reflections. The space of orbits $V/W$ of a finite Coxeter group is an affine variety; its coordinate ring coincides with the ring of $W$-invariant polynomials on $V$ \cite{Che}. The generators of such a ring are not uniquely specified, but their degrees $d_1=2<d_2\leq\ldots\leq d_{n-1}<d_n=h$ are invariants of the group;\footnote{In fact all inequalities are strict for all Coxeter groups but $D_n$.} the highest degree $h$ is called Coxeter number. \\
The complexification of the orbit space $M=(V \otimes \mathbb{C})/W$ can be endowed with a structure of Frobenius manifold. The two marked flat metrics have been worked out by K.Saito et al.\cite{Sa79,Sa93,SYS}, while the full construction is due to B.Dubrovin \cite{Du96}. We refer to his article for proofs of the following statements. \\
Let $x=(x_1,..,x_n)$ be coordinates on $V$ and $p=(p_{d_1},\ldots,p_{d_n})$ be local coordinates on $M$, where $p_{d_i}=p_{d_i}(x)$ are homogeneous generators of the ring of $W$-invariant polynomials. The second order generator $p_2$ fixes a $W$-invariant Euclidean metric $G_{ij}$ on $V$ via $p_2(x)=\frac{1}{2}\sum_{i,j}G_{ij}x_{i}x_{j}$. This induces a contravariant flat metric on $T_*M$ \cite{Sa79,Du96} by:
\begin{equation}
\label{met}
g^{\alpha\beta}(p)={<dp_{d_{\alpha}},dp_{d_{\beta}}>}^*=\sum_{i,j=1}^n \frac{\partial p_{d_{\alpha}}}{\partial x_i}\frac{\partial p_{d_{\beta}}}{\partial x_j}G^{ij},
\end{equation}
where $(G^{ij})=(G_{ij})^{-1}$. On the other hand, $\eta^{\alpha\beta}$ is defined by:
$$\eta^{\alpha\beta}(p)=\frac{\partial g^{\alpha\beta}(p)}{\partial p_h}.$$
The matrix of $\eta^{\alpha\beta}$ becomes constant once we calculate (Saito) flat coordinates $(t_1,\ldots,t_n)$. Dubrovin showed the existence of a relation between $g^{\alpha\beta}(t)$ and the Frobenius potential \cite{Du96}, given by:
\begin{equation}\label{sec}
g^{\alpha\beta}(t)={<dt_{d_{\alpha}},dt_{d_{\beta}}>}^*=\frac{(d_{\alpha}+d_{\beta}-2)}{h}
\eta^{\alpha\lambda}\eta^{\beta\mu}\partial_{\lambda}\partial_{\mu}F(t).
\end{equation}
This formula will allow us to reconstruct the potential (and hence the full Frobenius structure) from the knowledge of flat coordinates.

\section{Generators of $W$-invariant polynomials}
Systems of generators of $W$-invariant polynomials are well known for all irreducible finite Coxeter groups but the exceptional ones, labelled by $E_6$, $E_7$, $E_8$, for which very few and often cumbersome examples exist in the literature. In this note we use a construction due to M.Mehta \cite{Meh}, that has the advantage of being very natural to work out. For any $W$, we find a set of linear forms that is $W$-invariant, in the sense that, under the action of the group on the algebra of polynomials, the elements of the set transform into each other, leaving unchanged the whole set. Then we consider symmetric functions on these forms in the required degrees.

\subsection{Basic polynomials for $E_6$}
The action of the group $E_6$ on $\mathbb{R}^8$ is generated by reflections in the six hyperplanes $x_1=x_2$, $x_2=x_3$,\ldots,$x_5=x_6$ and $x_1+x_2+x_3+x_7=x_4+x_5+x_6+x_8$, all operating on the six dimensional subspace $S_6=\sum_{i=1}^6 x_i=0$, $x_7+x_8=0$. We first look for a set of polynomials in 8 variables, symmetric in $x_1,\ldots,x_6$ and in $x_7, x_8$ and invariant under reflection in the hyperplane $x_1+x_2+x_3+x_7=x_4+x_5+x_6+x_8$.
Apart from the obvious choices of $S_6$ and $x_7+x_8$, the following linear forms satisfy the above requirements:
$$\pm\frac{1}{2}(x_7-x_8)-x_i+\frac{1}{6}S_6, \qquad \qquad 1 \leq i \leq 6$$
$$x_i+x_j-\frac{1}{3}S_6, \qquad \qquad \qquad 1\leq i<j \leq 6$$
So the polynomials of the form:
\begin{align*}
u_m&=\sum_{i=1}^6\left(\left(\frac{1}{2}(x_7-x_8)-x_i+\frac{1}{6}S_6\right)^m +\left(-\frac{1}{2}(x_7-x_8)-x_i+\frac{1}{6}S_6\right)^m\right)+ \\
&+\sum_{1\leq i<j\leq 6}\left(x_i+x_j-\frac{1}{3}S_6\right)^m
\end{align*}
are manifestly invariant for any $m$. Restricting these polynomials to $S_6=0$, $x_7+x_8=0$ yields the standard realization of $E_6$ \cite{Meh}.
We recall that the $E_6$-invariant degrees are 2, 5, 6, 8, 9 and 12 \cite{Cox1}.
It turns out \cite{Cox2} that $u_2$, $u_5$, $u_6$, $u_8$, $u_9$ and $u_{12}$ are algebraically independent and form a system of generators for the ring of invariant polynomials of the $E_6$ type. We normalize the second order polynomial as $u_2/12$ (we will continue to call it $u_2$ for convenience), in such a way that it has as small as possible integer coefficients.

\subsection{Basic polynomials for $E_7$}
The action of the group $E_7$ on $\mathbb{R}^8$ is generated by reflections in the seven hyperplanes $x_1=x_2$, $x_2=x_3$,\ldots,$x_6=x_7$ and $x_1+x_2+x_3+x_4=x_5+x_6+x_7+x_8$, all operating on the seven dimensional subspace $S_7=\sum_{i=1}^8 x_i=0$. We first look for a set of polynomials in 8 variables, symmetric in $x_1,\ldots,x_8$ and invariant under reflection in the hyperplane $x_1+x_2+x_3+x_4=x_5+x_6+x_7+x_8$. Apart from the obvious choice of $S_7$, the following linear forms satisfy the above requirements:
$$\pm(x_i+x_j)-\frac{1}{4} S_7, \qquad 1 \leq i<j \leq 8$$
So the polynomials of the form:
\begin{displaymath}
v_m=\frac{1}{2}\sum_{1\leq i<j\leq 8}\left(\left(x_i+x_j-\frac{1}{4}S_7\right)^m+\left(-x_i-x_j-\frac{1}{4}S_7\right)^m\right)
\end{displaymath}
are manifestly invariant for any $m$. Restricting the action of the group on $S_7=0$ yields the standard realization of $E_7$ \cite{Meh}. We recall that the $E_7$-invariant degrees are 2, 6, 8, 10, 12, 14 and 18 \cite{Cox1}. It turns out \cite{Meh} that $v_2$, $v_6$, $v_8$, $v_{10}$, $v_{12}$, $v_{14}$ and $v_{18}$ are algebraically independent and form a system of generators for the ring of invariant polynomials of the $E_7$ type. We normalize the second order polynomial as $v_2/60$ (we will continue to call it $v_2$ for convenience), in such a way that it has as small as possible integer coefficients.

\subsection{Basic polynomials for $E_8$}
The action of the group $E_8$ on $\mathbb{R}^9$ is generated by reflections in the eight hyperplanes $x_1=x_2$, $x_2=x_3$,\ldots,$x_7=x_8$ and $2x_1+2x_2+2x_3=x_4+x_5+\ldots+x_9$, all operating on the eight dimensional subspace $S_8=\sum_{i=1}^9 x_i=0$. We first look for a set of polynomials in 9 variables, symmetric in $x_1,\ldots,x_8$ and invariant under reflection in the hyperplane $2x_1+2x_2+2x_3=x_4+x_5+\ldots+x_9$. Apart from the obvious choice of $S_8$, the following linear forms satisfy the above requirements:
$$\pm(x_i-x_j), \qquad \qquad \qquad \qquad 1 \leq i<j \leq 9$$
$$\pm\left(\frac{1}{3}S_8-x_i-x_j-x_k\right), \qquad \qquad 1 \leq i<j<k \leq 9$$
So the polynomials of the form:
\begin{displaymath}
w_m=\sum_{i<j}\left(x_i-x_j\right)^m+\sum_{i<j<k}\left(\frac{1}{3}S_8-x_i-x_j-x_k\right)^m
\end{displaymath}
are manifestly invariant for any even $m$. Restricting the action of the group on $S_8=0$ yields the standard realization of $E_8$ \cite{Meh}. We recall that the $E_8$-invariant degrees are 2, 8, 12, 14, 18, 20, 24 and 30 \cite{Cox1}. It turns out \cite{Meh} that $w_2$, $w_8$, $w_{12}$, $w_{14}$, $w_{18}$, $w_{20}$, $w_{24}$ and $w_{30}$ are algebraically independent and form a system of generators for ring of invariant polynomials of the $E_8$ type. We normalize the second order polynomial as $w_2/30$ (we will continue to call it $w_2$ for convenience), in such a way that it has as small as possible integer coefficients.

\section{Flat coordinates}
We present the explicit calculations needed to find the Frobenius structure for $E_6$. Analogous procedures will allow us to obtain the structures for $E_7$ and $E_8$. \\
We first need to write the metric (\ref{met}) with respect to the generators proposed by Mehta. For each matrix element ${<du_{d_{\alpha}},du_{d_{\beta}}>}^*$, we start finding all possible monic monomials in $(u_2,\ldots,u_{12})$ in the needed degree. Any such monomial $b$ must be written in the form:
\begin{displaymath}
b_{(i,c)}=\sum_{\substack{c_1 d_{i_1}+\ldots+c_s d_{i_s} \\ =d_{\alpha}+d_{\beta}-2}} u_{d_{i_1}}^{c_1}\ldots u_{d_{i_s}}^{c_s},
\end{displaymath}
where $i=(i_1,\ldots,i_s)$ and $c=(c_1,\ldots,c_s)$.\footnote{Note that in general this set of parameters uniquely identifies $b_{(i,c)}$ for all Coxeter groups but $D_n$.} \\
The next step consists in finding the (rational) coefficients $a_{(c,i)}$'s of the linear combination:
\begin{equation}\label{go}
{<du_{d_{\alpha}},du_{d_{\beta}}>}^*=\sum_{(c,i)}a_{(c,i)}b_{(c,i)},
\end{equation}
where the sum is over all allowed sets of parameters $(c,i)$.
This calculation is crucial from a computational point of view in our work; in fact, due to the high degrees of the polynomials and to the number of variables involved, standard computational software programs do not seem to be able to manage easily such an amount of information, that increases in complexity much more than linearly with the degree of the involved polynomials.\footnote{For example, ${<dw_{30},dw_{30}>}^*$ is a 58th degree homogeneous polynomial in 8 variables that has to be written as a linear combination of 163 polynomials of the same kind.} For this reason, we evaluate the expressions (\ref{go}) at points $(x_1,\ldots,x_n)\in \mathbb{N}^n$ such that $1\leq x_1\leq \ldots\leq x_n\leq n-2$. It turns out that, for $k$ generic choices of the $n$-tuple $(x_1,\ldots,x_n)$,\footnote{For example, $k=163$ for ${<dw_{30},dw_{30}>}^*$.} this procedure gives a - much easier to manage - determined linear system of equations in the $a_{(c,i)}$'s.
In Appendix A we present the matrix elements of the metric $g^{\alpha\beta}$ with respect to the Mehta generators for $E_6$, from which we can calculate the elements of the matrix $\eta^{\alpha\beta}$ different from zero:
\begin{align*}
&\partial_{u_{12}}{<du_2,du_{12}>}^*=\medmath{24} \qquad \qquad \partial_{u_{12}}{<du_5,du_9>}^*=\medmath{168} \\
&\partial_{u_{12}}{<du_6,du_8>}^*=\medmath{128} \qquad \qquad \partial_{u_{12}}{<du_6,du_{12}>}^*=\medmath{2752}\, u_2^2 \\
&\partial_{u_{12}}{<du_8,du_8>}^*=\medmath{896}\, u_2 \qquad \qquad
\partial_{u_{12}}{<du_8,du_{12}>}^*={\textstyle\frac{1064}{9}} \,u_6+{\textstyle\frac{25376}{3}}\, u_2^3 \\
&\partial_{u_{12}}{<du_9,du_9>}^*={\textstyle\frac{14112}{5}}\, u_2^2 \qquad \qquad \partial_{u_{12}}{<du_9,du_{12}>}^*=\medmath{742}\, u_5 u_2 \\
&\partial_{u_{12}}{<du_{12},du_{12}>}^*=\medmath{1254}\, u_8 u_2-{\textstyle\frac{6952}{9}}\, u_6 u_2^2+{\textstyle\frac{242}{5}}\, u_5^2+{\textstyle\frac{183656}{3}}\, u_2^5
\end{align*}
To find a flat basis, we now write the general homogeneous polynomials of degrees 2, 5, 6, 8,
9 and 12 belonging to our coordinate ring:
\begin{align*}
&t_2=u_2 \qquad t_5=u_5 \qquad t_6=u_6+k_1\, u_2^3 \qquad t_8=u_8+k_2\, u_6 u_2+k_3\, u_2^4 \\
&t_9=u_9+k_4\, u_5 u_2^2 \qquad t_{12}=u_{12}+k_5\, u_8 u_2^2+k_6\, u_6^2+k_7\, u_6 u_2^3+k_8\, u_5^2 u_2+k_9\, u_2^6
\end{align*}
where $k_i$ are free coefficients. Calculating the matrix elements of $\eta^{\alpha\beta}$ with respect to these new coordinates and asking for all the terms but the antidiagonal ones to be zero, we find an overdetermined system of equations that provides the coefficients needed to construct the Saito coordinates. We prefer to normalize some of the flat polynomials in such a way that all the elements in the antidiagonal are equal, as usual in the literature.
\begin{align*}
&t_2=u_2 \qquad \qquad t_5=u_5 \qquad \qquad t_6=u_6-\medmath{15}\,u_2^3 \\
&t_8={\textstyle\frac{3}{16}}(u_8-{\textstyle\frac{7}{2}}\, u_6 u_2+{\textstyle\frac{385}{12}}\, u_2^4) \qquad \qquad t_9={\textstyle\frac{1}{7}}\,(u_9-{\textstyle\frac{42}{5}}\, u_5 u_2^2) \\
&t_{12}=u_{12}-{\textstyle\frac{209}{16}}\, u_8 u_2^2-{\textstyle\frac{77}{576}}\, u_6^2+{\textstyle\frac{2959}{96}}\, u_6 u_2^3-{\textstyle\frac{121}{120}}\, u_5^2 u_2-{\textstyle\frac{6633}{32}}\, u_2^6 \\
\end{align*}
This result is already available in \cite{SYS} with different coefficients, because of a different choice of coordinates on $V$ done by the authors. \\
At this point we can calculate the matrix elements of $g^{\alpha\beta}(t)$ and, using \eqref{sec}, the Frobenius potential:
\begin{align*}
&F_{E_6}(t_2,t_5,t_6,t_8,t_9,t_{12}) ={\textstyle\frac{1}{24}}\,t_{12}^2t_2 +{\textstyle\frac{1}{12}}\,t_{12}t_8t_6 +{\textstyle\frac{1}{12}}\,t_{12}t_9t_5 +{\textstyle\frac{25}{147456}}\,t_6^4t_2+{\textstyle\frac{5}{3}}\,t_9^2t_8+ \\
&+{\textstyle\frac{25}{86016}}\,t_6^2t_2^7 +{\textstyle\frac{1}{38400}}\,t_6t_5^4 +{\textstyle\frac{1}{12800}}\,t_5^4t_2^3 +{\textstyle\frac{1}{8192}}\,t_5^2t_2^8 +{\textstyle\frac{1}{2048}}\,t_6t_5^2t_2^5 +{\textstyle\frac{1}{2048}}\,t_6^2t_5^2t_2^2+ \\
&+{\textstyle\frac{5}{768}}\,t_9t_6^2t_5 +{\textstyle\frac{1}{480}}\,t_9t_5^3t_2 +{\textstyle\frac{5}{384}}\,t_8t_6^2t_2^3 +{\textstyle\frac{5}{384}}\,t_9t_6t_5t_2^3 +{\textstyle\frac{1}{320}}\,t_8t_5^2t_2^4 +{\textstyle\frac{25}{192}}\,t_9^2t_2^4+ \\ &+{\textstyle\frac{1}{160}}\,t_8t_6t_5^2t_2 +{\textstyle\frac{25}{96}}\,t_9^2t_6t_2 +{\textstyle\frac{1}{50}}\,t_8^2t_5^2 +{\textstyle\frac{1}{20}}\,t_8^2t_2^5 +{\textstyle\frac{8}{15}}\,t_8^3t_2 +{\textstyle\frac{1}{4}}\,t_9t_8t_5t_2^2 +{\textstyle\frac{25}{1757184}}\,t_2^{13}
\end{align*}

\subsection{Frobenius structure for $E_7$}
Starting from the system of generators proposed by Mehta for $E_7$, we obtain:
\begin{align*}
t_2&=v_2 \\
t_6&=v_6-{\textstyle\frac{140}{9}}\,v_2^3 \\
t_8&=v_8-{\textstyle\frac{112}{27}}\,v_6 v_2+{\textstyle\frac{10220}{243}}\,v_2^4 \\
t_{10}&={\textstyle\frac{1}{70}}\,\big(v_{10}-{\textstyle\frac{9}{2}}\,v_8 v_2+\medmath{7}\, v_6 v_2^2-\medmath{42}\, v_2^5\big) \\
t_{12}&={\textstyle\frac{1}{1800}}\,\big(v_{12}-{\textstyle\frac{121}{21}}\, v_{10} v_2+{\textstyle\frac{341}{28}}\, v_8 v_2^2-{\textstyle\frac{11}{48}}\, v_6^2-{\textstyle\frac{649}{162}}\, v_6 v_2^3-{\textstyle\frac{8998}{729}}\, v_2^6\big) \\
t_{14}&={\textstyle\frac{1}{4466}}\,\big(v_{14}-{\textstyle\frac{7826}{1215}}\, v_{12} v_2+{\textstyle\frac{119977}{7290}}\, v_{10} v_2^2-{\textstyle\frac{1001}{2592}}\, v_8 v_6-{\textstyle\frac{434863}{29160}}\, v_8 v_2^3+ \\
&+{\textstyle\frac{253253}{174960}}\,v_6^2 v_2-{\textstyle\frac{18187169}{787320}}\,v_6 v_2^4+{\textstyle\frac{62327551}{354294}}\,v_2^7\big) \\
t_{18}&={\textstyle\frac{2}{1229}}\,\big(v_{18} -{\textstyle\frac{31144}{957}}\,v_{14}v_2^2 -{\textstyle\frac{2363}{4860}}\, v_{12}v_6+{\textstyle\frac{479381827}{3488265}}\,v_{12}v_2^3-{\textstyle\frac{3179}{8400}}\,v_{10}v_8+ \\
&+{\textstyle\frac{71893}{510300}}\,v_{10}v_6 v_2-{\textstyle\frac{1176266125}{5327532}}\,v_{10} v_2^4+{\textstyle\frac{14671}{50400}}\,v_8^2 v_2+{\textstyle\frac{32304709}{2466450}}\,v_8 v_6 v_2^2+ \\
&+{\textstyle\frac{43469918}{11099025}}\,v_8 v_2^5+{\textstyle\frac{117827}{2099520}}\,v_6^3 -{\textstyle\frac{272167739}{8456400}}\,v_6^2 v_2^3+{\textstyle\frac{301885587359}{513726300}}\,v_6 v_2^6- \\
&-{\textstyle\frac{4178016043387}{1387061010}}\,v_2^9\big) \\
\end{align*}

\newpage
\noindent
The resulting Frobenius potential is:
\begin{align*}
&F_{E_7}(t_2,t_6,t_8,t_{10},t_{12},t_{14},t_{18})= {\textstyle\frac{1}{36}}\,t_{18}^2t_2 +{\textstyle\frac{1}{36}}\,t_{18}t_{10}^2 +{\textstyle\frac{1}{18}}\,t_{18}t_{14}t_6 +{\textstyle\frac{1}{18}}\,t_{18}t_{12}t_8+ \\
&+{\textstyle\frac{800}{3}}\,t_{14}t_{12}^2 +\medmath{12}\,t_{14}^2t_{10} -{\textstyle\frac{16000}{243}}\,t_{12}^3 t_2 +{\textstyle\frac{80}{9}}\,t_{14}t_{12}t_{10}t_2  +{\textstyle\frac{400}{243}}\,t_{12}^2t_{10}t_2^2 +{\textstyle\frac{2}{27}}\,t_{14}^2t_6t_2^2 \\
&+{\textstyle\frac{2}{27}}\,t_{14}t_{10}^2t_2^2 +{\textstyle\frac{1}{6561}}\,t_{10}^3t_2^4 +{\textstyle\frac{16}{405}}\,t_{14}^2t_2^5 +{\textstyle\frac{1600}{413343}}\,t_{12}^2t_2^7+{\textstyle\frac{2}{4782969}}\,t_{10}^2t_2^9 +{\textstyle\frac{5}{243}}\,t_{12}t_{10}^2t_6 + \\
&+{\textstyle\frac{1}{4374}}\,t_{10}^3t_6t_2 +{\textstyle\frac{40}{729}}\,t_{14}t_{12}t_6t_2^3 -{\textstyle\frac{100}{19683}}\,t_{12}^2t_6t_2^4 +{\textstyle\frac{1}{2187}}\,t_{14}t_{10}t_6t_2^4 +{\textstyle\frac{1}{233280}}\,t_{14}t_8t_6^2t_2^2+ \\
&+{\textstyle\frac{1}{708588}}\,t_{10}^2t_6t_2^6+{\textstyle\frac{5}{486}}\,t_{14}t_{12}t_6^2 +{\textstyle\frac{50}{6561}}\,t_{12}^2t_6^2t_2 +{\textstyle\frac{1}{1944}}\,t_{14}t_{10}t_6^2t_2 +{\textstyle\frac{5}{26244}}\,t_{12}t_{10}t_6^2t_2^2 + \\ &+{\textstyle\frac{1}{229582512}}\,t_{10}t_6^2t_2^8 +{\textstyle\frac{1}{1889568}}\,t_{10}^2t_6^3 +{\textstyle\frac{1}{2834352}}\,t_{14}t_6^3t_2^3 +{\textstyle\frac{1}{6480}}\,t_{10}^3t_8 -{\textstyle\frac{2}{27}}\,t_{14}t_{12}t_8t_2^2+ \\
&+{\textstyle\frac{1}{102036672}}\,t_{10}t_6^3t_2^5 +{\textstyle\frac{1}{111577100832}}\,t_6^3t_2^{10} +{\textstyle\frac{5}{68024448}}\,t_{12}t_6^4t_2 +{\textstyle\frac{7}{110199605760}}\,t_8t_6^3t_2^6+ \\
&+{\textstyle\frac{1}{144000}}\,t_{14}t_8^3 +{\textstyle\frac{7}{528958107648}}\,t_6^5t_2^4 +{\textstyle\frac{1}{528958107648}}\,t_6^6t_2 +{\textstyle\frac{1}{5}}\,t_{14}^2t_8t_2 +{\textstyle\frac{1}{1080}}\,t_{14}t_{10}t_8t_6+ \\
& +{\textstyle\frac{20}{729}}\,t_{12}^2t_8t_2^3 +{\textstyle\frac{1}{810}}\,t_{14}t_{10}t_8t_2^3 +{\textstyle\frac{1}{2187}}\,t_{12}t_{10}t_8t_2^4 -{\textstyle\frac{5}{486}}\,t_{12}^2t_8t_6 -{\textstyle\frac{1}{10203667200}}\,t_8^2t_6t_2^8+ \\
&+{\textstyle\frac{1}{2916}}\,t_{12}t_{10}t_8t_6t_2 +{\textstyle\frac{1}{116640}}\,t_{10}^2t_8t_6t_2^2 +{\textstyle\frac{1}{708588}}\,t_{12}t_8t_6t_2^6 +{\textstyle\frac{1}{42515280}}\,t_{10}t_8t_6t_2^7 + \\
&+{\textstyle\frac{1}{22674816}}\,t_{10}t_8t_6^2t_2^4 +{\textstyle\frac{1}{13774950720}}\,t_8t_6^2t_2^9 +{\textstyle\frac{1}{2519424}}\,t_{12}t_8t_6^3 +{\textstyle\frac{1}{50388480}}\,t_{10}t_8t_6^3t_2 + \\
&+{\textstyle\frac{1}{31492800}}\,t_{10}t_8^2t_2^6 +{\textstyle\frac{1}{21045063600}}\,t_8^2t_2^{11} +{\textstyle\frac{1}{116640}}\,t_{12}t_8^2t_6t_2^2 +{\textstyle\frac{1}{20995200}}\,t_{10}t_8^2t_6t_2^3 + \\
&+{\textstyle\frac{1}{55987200}}\,t_{10}t_8^2t_6^2 +{\textstyle\frac{1}{2267481600}}\,t_8^2t_6^2t_2^5 +{\textstyle\frac{1}{157464000}}\,t_{14}t_8^2t_6^2t_2^6 +{\textstyle\frac{1}{18139852800}}\,t_8^2t_6^3t_2^2+ \\
&+{\textstyle\frac{1}{5668704000}}\,t_8^3t_2^7 -{\textstyle\frac{1}{3023308800}}\,t_8^3t_6t_2^4 +{\textstyle\frac{1}{3359232000}}\,t_8^3t_6^2t_2 +{\textstyle\frac{1}{1399680000}}\,t_8^4t_2^3 - \\
&-{\textstyle\frac{1}{7464960000}}\,t_8^4t_6 -{\textstyle\frac{1}{194400}}\,t_{12}t_8^3t_2 +{\textstyle\frac{1}{181398528}}\,t_{10}t_6^4t_2^2+{\textstyle\frac{1}{181312788852}}\,t_6^2t_2^{13}+ \\
&+{\textstyle\frac{7}{297538935552}}\,t_6^4t_2^7 +{\textstyle\frac{1}{314928}}\,t_{10}^2t_6^2t_2^3 +{\textstyle\frac{25}{76527504}}\,t_{12}t_6^3t_2^4 +{\textstyle\frac{10}{59049}}\,t_{12}t_{10}t_6t_2^5+ \\ &+{\textstyle\frac{1}{10077696}}\,t_{14}t_6^4 +{\textstyle\frac{1}{194400}}\,t_{14}t_8^2t_2^4 -{\textstyle\frac{1}{262440}}\,t_{12}t_8^2t_2^5 +{\textstyle\frac{1}{261213880320}}\,t_8t_6^5 + \\
&+{\textstyle\frac{5}{58773123072}}\,t_8t_6^4t_2^3 +{\textstyle\frac{1}{129600}}\,t_{10}^2t_8^2t_2 +{\textstyle\frac{16}{91224740283363}}\,t_2^{19}
\end{align*}

\subsection{Frobenius structure for $E_8$}
Starting from the system of generators proposed by Mehta for $E_8$, we obtain:
\begin{align*}
t_2&=w_2 \\
t_8&=w_8-{\textstyle\frac{1176}{5}}\,w_2^4 \\
t_{12}&=w_{12}-{\textstyle\frac{363}{10}}\,w_8w_2^2+{\textstyle\frac{924924}{125}}\,w_2^6 \\
t_{14}&=w_{14}-{\textstyle\frac{169}{18}}\,w_{12}w_2+{\textstyle\frac{35321}{180}}\,w_8w_2^3-{\textstyle\frac{2859142}{75}}\, w_2^7 \\
t_{18}&={\textstyle\frac{5}{42}}\,\big(w_{18}-{\textstyle\frac{867}{11}}\,w_{14}w_2^2+{\textstyle\frac{161551}{330}}\,w_{12}w_2^3 -{\textstyle\frac{3757}{1920}}\,w_8^2w_2 -{\textstyle\frac{7172113}{1000}}\,w_8w_2^5+ \\
&+{\textstyle\frac{1789913697}{1250}}\,w_2^9\big) \\
\end{align*}
\begin{align*}
t_{20}&={\textstyle\frac{325}{2091}}\,\big(w_{20} -{\textstyle\frac{221293}{16380}}\,w_{18}w_2+{\textstyle\frac{509371}{780}}\,w_{14}w_2^3
-{\textstyle\frac{6137}{25920}}\,w_{12}w_8 -{\textstyle\frac{11666437}{3150}}\,w_{12}w_2^4+ \\ &+{\textstyle\frac{16134173}{806400}}\,w_8^2w_2^2 +{\textstyle\frac{12765371179}{252000}}\,w_8w_2^6-{\textstyle\frac{48701290603009}{4725000}}\,w_2^{10}\big) \\
t_{24}&={\textstyle\frac{1625}{15124}}\,\big(w_{24}-{\textstyle\frac{4558393}{34850}}\,w_{20}w_2^2 +{\textstyle\frac{278965917091}{244647000}}\,w_{18}w_2^3 -{\textstyle\frac{412091}{72000}}\,w_{14}w_8w_2- \\
&-{\textstyle\frac{2520034190873}{59962500}}\,w_{14}w_2^5-{\textstyle\frac{130663}{1360800}}\,w_{12}^2 +{\textstyle\frac{2147064467}{37195200}}\,w_{12}w_8w_2^2+ \\
&+{\textstyle\frac{37784188712363}{166050000}}\,w_{12}w_2^6-{\textstyle\frac{110561}{4976640}}\,w_8^3 -{\textstyle\frac{5445874541701}{2975616000}}\,w_8^2w_2^4- \\
&-{\textstyle\frac{102233301029629}{34593750}}\,w_8w_2^8 +{\textstyle\frac{288133137433526381}{461250000}}\,w_2^{12}\big) \\
t_{30}&={\textstyle\frac{96}{61}}\,\big(w_{30} -{\textstyle\frac{830066159}{665456}}\,w_{24}w_2^3-{\textstyle\frac{109362799}{14720640}}\,w_{20}w_8w_2 +{\textstyle\frac{3331116694332609}{23191141600}}\,w_{20}w_2^5- \\
&-{\textstyle\frac{15979}{72576}}\,w_{18}w_{12}+{\textstyle\frac{216507881}{5152224}}\,w_{18}w_8w_2^2 -{\textstyle\frac{264018066617724437}{219705552000}}\,w_{18}w_2^6- \\
&-{\textstyle\frac{9209791}{2710400}}\,w_{14}^2w_2 +{\textstyle\frac{429264061}{14636160}}\,w_{14}w_{12}w_2^2-{\textstyle\frac{10933}{138240}}\,w_{14}w_8^2+ \\
&+{\textstyle\frac{1250583443140331}{241245312000}}\,w_{14}w_8w_2^4+{\textstyle\frac{1084439719582395161}{25129720000}}\,w_{14}w_2^8+ \\
&+{\textstyle\frac{706407401999}{10485345024}}\,w_{12}^2w_2^3+{\textstyle\frac{535397745467}{314225049600}}\,w_{12}w_8^2w_2 -{\textstyle\frac{20688389046489203}{361867968000}}\,w_{12}w_8w_2^5- \\
&-{\textstyle\frac{2040455724082448777}{8767440000}}\,w_{12}w_2^9 -{\textstyle\frac{1521588819282337}{113692336128000}}\,w_8^3w_2^3 +{\textstyle\frac{9670511305095824287}{5263534080000}}\,w_8^2w_2^7+ \\
&+{\textstyle\frac{640920923508470286331}{211481280000}}\,w_8w_2^{11}-{\textstyle\frac{130055065986893806638453467}{203550732000000}}\,w_2^{15}\big) \\
\end{align*}

\noindent
The resulting Frobenius potential is:
\begin{align*}
&F_{E_8}(t_2,t_8,t_{12},t_{14},t_{18},t_{20},t_{24},t_{30})={\textstyle\frac{1}{60}}\,t_{30}^2t_2 +{\textstyle\frac{1}{30}}\,t_{30}t_{20}t_{12} +{\textstyle\frac{1}{30}}\,t_{30}t_{24}t_8+ \\
&+{\textstyle\frac{1}{30}}\,t_{30}t_{18}t_{14} +{\textstyle\frac{11}{540}}\,t_{18}^2t_{14}t_{12} +{\textstyle\frac{20449}{3888730944}}\, t_{12}^5t_2 +{\textstyle\frac{13}{6930}}\,t_{18}t_{14}^3t_2 +{\textstyle\frac{169}{873180}}\,t_{14}^3t_{12}t_2^4+ \\
&+{\textstyle\frac{121}{4860}}\,t_{18}^2t_{12}^2t_2+{\textstyle\frac{13}{1890}}\,t_{18}t_{14}^2t_{12}t_2^2 +{\textstyle\frac{169}{1067220}}\,t_{14}^4t_2^3+{\textstyle\frac{143}{25515}}\,t_{18}t_{14}t_{12}^2t_2^3 + \\
&+{\textstyle\frac{169}{396900}}\,t_{14}^2t_{12}^2t_2^5+{\textstyle\frac{1232}{16875}}\,t_{18}^2t_{14}t_2^6 +{\textstyle\frac{20449}{289340100}}\,t_{12}^4t_2^7+{\textstyle\frac{3388}{151875}}\,t_{18}^2t_{12}t_2^7 +{\textstyle\frac{1024}{4375}}\,t_{24}^2t_2^7 + \\
&+{\textstyle\frac{1144}{151875}}\,t_{18}t_{14}t_{12}t_2^9+{\textstyle\frac{2704}{19490625}}\,t_{14}^3t_2^{10} +{\textstyle\frac{676}{5315625}}\,t_{14}^2t_{12}t_2^{11}+{\textstyle\frac{81796}{1291696875}}\,t_{12}^3t_2^{13}+ \\
&+{\textstyle\frac{1517824}{82265625}}\,t_{18}^2t_2^{13}+{\textstyle\frac{43264}{537890625}}\,t_{14}^2t_2^{17} +{\textstyle\frac{5234944}{146084765625}}\,t_{12}^2t_2^{19} +{\textstyle\frac{262668550144}{180003021240234375}}\,t_2^{31}+ \\ &+{\textstyle\frac{13}{8470}}\,t_{14}^3+{\textstyle\frac{13}{2310}}\,t_{20}t_{14}^2t_{12}t_2 +{\textstyle\frac{143}{51030}}\,t_{20}t_{12}^3t_2^3+{\textstyle\frac{2156}{1625}}\,t_{20}t_{18}^2t_2^3 +{\textstyle\frac{56}{625}}\,t_{20}t_{18}t_{14}t_2^5+ \\
&+{\textstyle\frac{52}{6875}}\,t_{20}t_{14}^2t_2^7+{\textstyle\frac{572}{151875}}\,t_{20}t_{12}^2t_2^9 +{\textstyle\frac{84}{625}}\,t_{20}^2t_{12}t_2^5+{\textstyle\frac{37632}{859375}}\,t_{20}^2t_2^{11} +{\textstyle\frac{15876}{17875}}\,t_{20}^3t_2+ \\
&+{\textstyle\frac{13}{3969}}\,t_{24}t_{14}t_{12}^2+{\textstyle\frac{616}{325}}\,t_{24}t_{18}^2t_2 +{\textstyle\frac{16}{75}}\,t_{24}t_{18}t_{14}t_2^3+{\textstyle\frac{176}{675}}\,t_{24}t_{18}t_{12}t_2^4 +{\textstyle\frac{104}{9625}}\,t_{24}t_{14}^2t_2^5+ \\
&+{\textstyle\frac{208}{23625}}\,t_{24}t_{14}t_{12}t_2^6+{\textstyle\frac{1008}{325}}\,t_{24}t_{20}t_{18} +{\textstyle\frac{144}{275}}\,t_{24}t_{20}t_{14}t_2^2+{\textstyle\frac{96}{385}}\,t_{24}^2t_{14} +{\textstyle\frac{32}{105}}\,t_{24}^2t_{12}t_2+ \\
& +{\textstyle\frac{1573}{8817984}}\,t_{18}t_{12}^3t_8 +{\textstyle\frac{169}{4191264}}\,t_{14}^3t_{12}t_8 +{\textstyle\frac{1859}{564350976}}\,t_{14}t_{12}^2t_8^3 +{\textstyle\frac{49}{2000}}\,t_{20}^2t_8^2t_2^3+ \\
&+{\textstyle\frac{5929}{70200}}\,t_{18}^3t_8+{\textstyle\frac{169}{2286144}}\,t_{14}^2t_{12}^2t_8t_2 +{\textstyle\frac{1859}{20575296}}\,t_{14}t_{12}^3t_8t_2^2+{\textstyle\frac{77}{1800}}\,t_{18}^2t_{14}t_8t_2^2+ \\
&+{\textstyle\frac{847}{24300}}\,t_{18}^2t_{12}t_8t_2^3+{\textstyle\frac{13}{1800}}\,t_{18}t_{14}^2t_8t_2^4 +{\textstyle\frac{143}{24300}}\,t_{18}t_{14}t_{12}t_8t_2^5+{\textstyle\frac{169}{891000}}\,t_{14}^3t_8t_2^6+ \\
&+{\textstyle\frac{11011}{3280500}}\,t_{18}t_{12}^2t_8t_2^6+{\textstyle\frac{169}{425250}}\,t_{14}^2t_{12}t_8t_2^7 +{\textstyle\frac{1859}{5103000}}\,t_{14}t_{12}^2t_8t_2^8+{\textstyle\frac{11858}{759375}}\,t_{18}^2t_8t_2^9+ \\
&+{\textstyle\frac{4004}{1265625}}\,t_{18}t_{14}t_8t_2^{11}+{\textstyle\frac{44044}{34171875}}\,t_{18}t_{12}t_8t_2^{12} +{\textstyle\frac{338}{3796875}}\,t_{14}^2t_8t_2^{13}+{\textstyle\frac{7436}{34171875}}\,t_{14}t_{12}t_8t_2^{14}+ \\
\end{align*}
\begin{align*}
&+{\textstyle\frac{128256128}{8649755859375}}\,t_8t_2^{27}+{\textstyle\frac{7}{100}}\,t_{20}t_{18}t_{14}t_8t_2 +{\textstyle\frac{77}{900}}\,t_{20}t_{18}t_{12}t_8t_2^2+{\textstyle\frac{13}{3300}}\,t_{20}t_{14}^2t_8t_2^3+ \\
&+{\textstyle\frac{13}{1350}}\,t_{20}t_{14}t_{12}t_8t_2^4+{\textstyle\frac{2156}{28125}}\,t_{20}t_{18}t_8t_2^8 +{\textstyle\frac{364}{140625}}\,t_{20}t_{14}t_8t_2^{10}+{\textstyle\frac{63}{2200}}\,t_{20}^2t_{14}t_8+ \\
&+{\textstyle\frac{11}{270}}\,t_{24}t_{18}t_{12}t_8+{\textstyle\frac{13}{2310}}\,t_{24}t_{14}^2t_8t_2 +{\textstyle\frac{13}{1890}}\,t_{24}t_{14}t_{12}t_8t_2^2+{\textstyle\frac{143}{51030}}\,t_{24}t_{12}^2t_8t_2^3+ \\
&+{\textstyle\frac{1232}{16875}}\,t_{24}t_{18}t_8t_2^6+{\textstyle\frac{112}{625}}\,t_{24}t_{20}t_8t_2^5 +{\textstyle\frac{16}{75}}\,t_{24}^2t_8t_2^3 +{\textstyle\frac{13}{34560}}\,t_{18}t_{14}^2t_8^2+ \\
&+{\textstyle\frac{143}{103680}}\,t_{18}t_{14}t_{12}t_8^2t_2 +{\textstyle\frac{169}{2661120}}\,t_{14}^3t_8^2t_2^2+{\textstyle\frac{1573}{2799360}}\,t_{18}t_{12}^2t_8^2t_2^2 +{\textstyle\frac{169}{933120}}\,t_{14}^2t_{12}t_8^2t_2^3+ \\
&+{\textstyle\frac{18592}{14696640}}\,t_{14}t_{12}^2t_8^2t_2^4+{\textstyle\frac{20449}{377913600}}\,t_{12}^3t_8^2t_2^5 +{\textstyle\frac{5929}{270000}}\,t_{18}^2t_8^2t_2^5+{\textstyle\frac{1001}{303750}}\,t_{18}t_{14}t_8^2t_2^7+ \\
&+{\textstyle\frac{11011}{7290000}}\,t_{18}t_{12}t_8^2t_2^8+{\textstyle\frac{169}{810000}}\,t_{14}^2t_8^2t_2^9 +{\textstyle\frac{20449}{109350000}}\,t_{14}t_{12}t_8^2t_2^{10}+{\textstyle\frac{20449}{196830000}}\,t_{12}^2t_8^2t_2^{11}+ \\
&+{\textstyle\frac{26026}{854296875}}\,t_{14}t_8^2t_2^{16}+{\textstyle\frac{286286}{7688671875}}\,t_{12}t_8^2t_2^{17} +{\textstyle\frac{64128064}{7368310546875}}\,t_8^2t_2^{23}+{\textstyle\frac{13}{34560}}\,t_{20}t_{14}t_{12}t_8^2+ \\
&+{\textstyle\frac{143}{103680}}\,t_{20}t_{12}^2t_8^2t_2+{\textstyle\frac{539}{18000}}\,t_{20}t_{18}t_8^2t_2^4 +{\textstyle\frac{637}{270000}}\,t_{20}t_{14}t_8^2t_2^6+{\textstyle\frac{1001}{607500}}\,t_{20}t_{12}t_8^2t_2^7 + \\
&+{\textstyle\frac{77}{1800}}\,t_{24}t_{18}t_8^2t_2^2+{\textstyle\frac{13}{5400}}\,t_{24}t_{14}t_8^2t_2^4 +{\textstyle\frac{143}{243000}}\,t_{24}t_{12}t_8^2t_2^5+{\textstyle\frac{7}{200}}\,t_{24}t_{20}t_8^2t_2 + \\
&+{\textstyle\frac{20449}{10158317568}}\,t_{12}^3t_8^3t_2+{\textstyle\frac{5929}{3110400}}\,t_{18}^2t_8^3t_2 +{\textstyle\frac{1001}{1166400}}\,t_{18}t_{14}t_8^3t_2^3+{\textstyle\frac{11011}{16796160}}\,t_{18}t_{12}t_8^3t_2^4+ \\
&+{\textstyle\frac{1183}{15552000}}\,t_{14}^2t_8^3t_2^5+{\textstyle\frac{13013}{104976000}}\,t_{14}t_{12}t_8^3t_2^6 +{\textstyle\frac{143143}{3779136000}}\,t_{12}^2t_8^3t_2^7+{\textstyle\frac{847847}{3280500000}}\,t_{18}t_8^3t_2^{10}+ \\
&+{\textstyle\frac{169169}{3280500000}}\,t_{14}t_8^3t_2^{12}+{\textstyle\frac{1002001}{29524500000}}\,t_{12}t_8^3t_2^{13} +{\textstyle\frac{1002001}{115330078125}}\,t_8^3t_2^{19}+{\textstyle\frac{539}{345600}}\,t_{20}t_{18}t_8^3+ \\
&+{\textstyle\frac{91}{172800}}\,t_{20}t_{14}t_8^3t_2^2+{\textstyle\frac{1001}{3110400}}\,t_{20}t_{12}t_8^3t_2^3 +{\textstyle\frac{7007}{72900000}}\,t_{20}t_8^3t_2^9+{\textstyle\frac{13}{103680}}\,t_{24}t_{14}t_8^3 + \\
&+{\textstyle\frac{11011}{1074954240}}\,t_{18}t_{12}t_8^4+{\textstyle\frac{169}{47775744}}\,t_{14}^2t_8^4t_2 +{\textstyle\frac{20449}{2149908480}}\,t_{14}t_{12}t_8^4t_2^2+{\textstyle\frac{347633}{58047528960}}\,t_{12}^2t_8^4t_2^3+ \\
&+{\textstyle\frac{539539}{4199040000}}\,t_{18}t_8^4t_2^6+{\textstyle\frac{13013}{466560000}}\,t_{14}t_8^4t_2^8 +{\textstyle\frac{11022011}{453496320000}}\,t_{12}t_8^4t_2^9+{\textstyle\frac{1002001}{92264062500}}\,t_8^4t_2^{15}+ \\
&+{\textstyle\frac{49049}{622080000}}\,t_{20}t_8^4t_2^5+{\textstyle\frac{1001}{37324800}}\,t_{24}t_8^4t_2^3 +{\textstyle\frac{77077}{7166361600}}\,t_{18}t_8^5t_2^2+{\textstyle\frac{91091}{21499084800}}\,t_{14}t_8^5t_2^4+ \\
&+{\textstyle\frac{1002001}{241864704000}}\,t_{12}t_8^5t_2^5+{\textstyle\frac{13026013}{2519424000000}}\,t_8^5t_2^{11} +{\textstyle\frac{7007}{1990656000}}\,t_{20}t_8^5t_2+{\textstyle\frac{13013}{515978035200}}\,t_{14}t_8^6+ \\
&+{\textstyle\frac{143}{933120}}\,t_{24}t_{12}t_8^3t_2 +{\textstyle\frac{143143}{1857520926720}}\,t_{12}t_8^6t_2+{\textstyle\frac{11022011}{9674588160000}}\,t_8^6t_2^7 +{\textstyle\frac{1002001}{18575209267200}}\,t_8^7t_2^3
\end{align*}

\appendix
\section{Matrix elements of $g^{\alpha\beta}$ for $E_6$}
In this Appendix we present the explicit calculations of $g^{\alpha\beta}(u_2,\ldots,u_{12})$ for $E_6$ with respect to the coordinates corresponding to Mehta polynomials.
\begin{align*}
{<du_2,du_k>}^*&=\medmath{2k}\, u_k \\
{<du_5,du_5>}^*&=\medmath{120}\, u_8-\medmath{320}\, u_6 u_2+\medmath{2400}\, u_2^4 \\
{<du_5,du_6>}^*&={\textstyle\frac{720}{7}}\,u_9-\medmath{360}\, u_5 u_2^2 \\
{<du_5,du_8>}^*&={\textstyle\frac{9760}{21}}\,u_9 u_2+{\textstyle\frac{160}{3}}\,u_6 u_5-\medmath{3680}\, u_5 u_2^3 \\
{<du_5,du_9>}^*&=\medmath{168}\, u_{12}-\medmath{1092}\, u_8 u_2^2-{\textstyle\frac{28}{3}}\,u_6^2+\medmath{1792}\, u_6 u_2^3+{\textstyle\frac{56}{5}}\,u_5^2 u_2^2-\medmath{9072}\, u_2^6 \\
\end{align*}
\begin{align*}
{<du_5,du_{12}>}^*&={\textstyle\frac{380}{7}}\,u_9 u_6+{\textstyle\frac{53600}{21}}\,u_9 u_2^3+\medmath{537}\, u_8 u_5 u_2- \\
&-{\textstyle\frac{3896}{3}}\,u_6 u_5 u_2^2+{\textstyle\frac{65}{6}}\,u_5^3-\medmath{15900}\, u_5 u_2^5 \\
{<du_6,du_6>}^*&=\medmath{144}\, u_8 u_2+\medmath{576}\, u_6 u_2^2+{\textstyle\frac{48}{5}}\,u_5^2-\medmath{3456}\, u_2^5\\
{<du_6,du_8>}^*&=\medmath{128}\, u_{12}-\medmath{448}\, u_8 u_2^2+{\textstyle\frac{224}{9}}\,u_6^2+ \\
&+{\textstyle\frac{5824}{3}}\,u_6 u_2^3-{\textstyle\frac{1288}{15}}\,u_5^2 u_2^2-\medmath{16128}\, u_2^6\\
{<du_6,du_9>}^*&=\medmath{1752}\, u_9 u_2^2+{\textstyle\frac{819}{25}}\,u_8 u_5+{\textstyle\frac{2856}{25}}\,u_6 u_5 u_2-\medmath{9828}\, u_5 u_2^4 \\
{<du_6,du_{12}>}^*&=\medmath{2752}\, u_{12}u_2^2+{\textstyle\frac{1644}{7}}\,u_9 u_5 u_2+{\textstyle\frac{108}{5}}\, u_8^2+{\textstyle\frac{1004}{5}}\,u_8 u_6 u_2- \\
&-\medmath{22976}\, u_8 u_2^4-{\textstyle\frac{22288}{45}}\,u_6^2 u_2^2+{\textstyle\frac{256}{15}}\,u_6 u_5^2+ \\
&+{\textstyle\frac{148400}{3}}\,u_6 u_2^5-{\textstyle\frac{52412}{15}}\,u_5^2 u_2^3-\medmath{301248}\, u_2^8 \\
{<du_8,du_8>}^*&=\medmath{896}\, u_{12}u_2+{\textstyle\frac{160}{9}}\,u_9 u_5+\medmath{112}\, u_8 u_6-\medmath{8960}\, u_8 u_2^3-{\textstyle\frac{448}{3}}\, u_6^2 u_2+ \\
&+\medmath{20160}\, u_6 u_2^4-{\textstyle\frac{2576}{3}}\,u_5^2 u_2^2-\medmath{130560}\, u_2^7 \\
{<du_8,du_9>}^*&=\medmath{88}\, u_9 u_6+\medmath{4352}\, u_9 u_2^3+{\textstyle\frac{10374}{25}}\,u_8 u_5 u_2- \\
&-{\textstyle\frac{17024}{25}}\,u_6 u_5 u_2^2+{\textstyle\frac{7}{5}}\,u_5^3-\medmath{31752}\, u_5 u_2^5 \\
{<du_8,du_{12}>}^*&=u_{12}\big( {\textstyle\frac{1064}{9}}\,u_6+{\textstyle\frac{25376}{3}}\,u_2^3\big)+{\textstyle\frac{600}{49}}\,u_9^2
+{\textstyle\frac{67576}{63}}\,u_9 u_5 u_2^2+ \\
&+{\textstyle\frac{2468}{5}}\,u_8^2 u^2-{\textstyle\frac{77804}{45}}\,u_8 u_6 u_2^2+{\textstyle\frac{113}{6}}\,u_8 u_5^2-{\textstyle\frac{264016}{3}}\,u_8u_2^5- \\
&-{\textstyle\frac{244}{81}}\,u_6^3+{\textstyle\frac{16304}{15}}\,u_6^2 u_2^3+{\textstyle\frac{452}{135}}\,u_6 u_5^2 u_2+{\textstyle\frac{1694512}{9}}\,u_6 u_2^6- \\
&-{\textstyle\frac{757426}{45}}\,u_5^2 u_2^4-\medmath{1281408}\, u_2^9 \\
{<du_9,du_9>}^*&={\textstyle\frac{14112}{5}}\,u_{12}u_2^2+{\textstyle\frac{3192}{5}}\,u_9 u_5 u_2+{\textstyle\frac{1323}{50}}\,u_8^2-{\textstyle\frac{3528}{25}}\,u_8 u_6 u_2- \\
&-{\textstyle\frac{128772}{5}}\,u_8 u_2^4+{\textstyle\frac{784}{25}}\,u_6^2 u_2^2+{\textstyle\frac{147}{25}}\,u_6 u_5^2+{\textstyle\frac{249312}{5}}\,u_6 u_2^5- \\
&-\medmath{4116}\, u_5^2 u_2^3-{\textstyle\frac{1555848}{5}}\,u_2^8 \\
{<du_9,du_{12}>}^*&=\medmath{742}\, u_{12}u_5 u_2 +\medmath{513}\, u_9 u_8 u_2+{\textstyle\frac{140}{3}}\,u_9 u_6 u_2^2+\medmath{21}\, u_9 u_5^2+ \\
&+\medmath{23412}\, u_9 u_2^5+{\textstyle\frac{406}{25}}\,u_8 u_6 u_5-{\textstyle\frac{16156}{5}}\,u_8 u_5 u_2^3-{\textstyle\frac{13979}{225}}\,u_6^2 u_5 u_2- \\
&-{\textstyle\frac{3752}{15}}\,u_6 u_5 u_2^4-{\textstyle\frac{1141}{3}}\,u_5^3 u_2^2 -\medmath{193032}\, u_5 u_2^7 \\
{<du_{12},du_{12}>}^*&=u_{12}\big( \medmath{1254}\, u_8 u_2-{\textstyle\frac{6952}{9}}\,u_6 u_2^2+{\textstyle\frac{242}{5}}\, u_5^2+{\textstyle\frac{183656}{3}}\,u_2^5 \big)+ \\
&+{\textstyle\frac{16500}{49}}\,u_9^2 u_2^2+{\textstyle\frac{99}{14}}\,u_9 u_8 u_5+{\textstyle\frac{4510}{21}}\,u_9 u_6 u_5 u_2+{\textstyle\frac{188210}{21}}\,u_9 u_5 u_2^4+ \\
&+{\textstyle\frac{231}{20}}\,u_8^2 u_6-\medmath{2068}\, u_8^2 u_2^3-{\textstyle\frac{2827}{45}}\,u_8 u_6^2 u_2-{\textstyle\frac{58498}{9}}\,u_8 u_6 u_2^4+ \\
&+{\textstyle\frac{6424}{15}}\,u_8 u_5^2 u_2^2-{\textstyle\frac{1885048}{3}}\,u_8 u_2^7+{\textstyle\frac{3124}{405}}\,u_6^3 u_2^2+{\textstyle\frac{671}{270}}\,u_6^2 u_5^2+ \\
&+{\textstyle\frac{53372}{3}}\,u_6^2 u_2^5-{\textstyle\frac{83072}{27}}\,u_6 u_5^2 u_2^3+{\textstyle\frac{10460428}{9}}\, u_6 u_2^8-{\textstyle\frac{143}{180}}\,u_5^4 u_2- \\
&-{\textstyle\frac{6338332}{45}}\,u_5^2 u_2^6-\medmath{8443104}\, u_2^{11}\\
\end{align*}


\begin{thebibliography}{99}

\bibitem[Ch55]{Che} Chevalley C., Invariants of finite groups generated by reflections, \emph{Amer. J. Math.} {\bf{77}} (1955), 778-782.

\bibitem[Co34]{Cox1} Coxeter H.S.M., Discrete groups generated by reflections, \emph{Ann. Math.} {\bf{35}} (1934), 588-621.

\bibitem[Co51]{Cox2} Coxeter H.S.M., The product of the generators of a finite group generated by reflections, \emph{Duke Math. J.} {\bf{18}} (1951), 765-782.

 \bibitem[DLZ93]{DLZ} Di Francesco P., Lesage F., Zuber J.-B., Graph rings and integrable perturbations of N=2 superconformal theories, \emph{Nucl. Phys.} {\bf{B408}} (1993), 600-634.

\bibitem[Du96]{Du96} Dubrovin B., Geometry of 2D topological field theories, \emph{Lect. Notes in Math.} {\bf{1620}} (1996), 120-348.

\bibitem[KW81]{Kat} Kato M. and Watanabe S., The flat coordinate system of the rational double point
of $E_8$ type, \emph{Bull. Coll. Sci., Univ. Ryukyus} {\bf{32}} (1981), 1-3.

\bibitem [Me88]{Meh} Mehta M.L., Basic set of invariant polynomials for finite reflection groups, \emph{Comm. Alg.}  {\bf{16}} (1988), 1083-1098.

\bibitem[Sa79]{Sa79} Saito K., Extended affine root systems II (flat invariants), \emph{Publ. RIMS, Kyoto Univ.} {\bf{26}} (1979), 15-78.

\bibitem [Sa93]{Sa93} Saito K., On a linear structure of the quotient variety by a finite reflexion group, \emph{Publ. RIMS, Kyoto Univ.} {\bf{29(4)}} (1993), 535-579.

\bibitem [SYS80]{SYS} Saito K., Yano T. and Sekiguchi J., On a certain generator system of the ring of invariants of a finite refection group, \emph{Comm. Alg.} {\bf{8}} (1980), 373-408.

\bibitem [Ya80] {Ya}Yano T., Free deformation for isolated singularity, \emph{Sci. Rep. Saitama Univ.} {\bf{A9}} (1980), 61-70.

\end{thebibliography}
\end{document}